\documentclass[12pt]{article}
\usepackage{epsfig,amsmath,amssymb,graphicx}
\usepackage{comment,fullpage}
\usepackage{authblk}

\newtheorem{lemma}{Lemma}[section]
\newtheorem{theorem}[lemma]{Theorem}

\newtheorem{conjecture}[lemma]{Conjecture}

\newtheorem{example}[lemma]{Example}

\newcommand {\qed}{\rule{2mm}{2mm}\medskip}
\newenvironment {proof}{\noindent{\bf Proof:}~}{~\qed}

\newcommand\Taxi{\operatorname{Taxicab}}

\begin{document}

% THE TITLE PAGE 

\title{ Seeds for Generalized Taxicab Numbers}

\author[1]{Jeffrey H. Dinitz}
\author[2]{Richard Games}
\author[3]{Robert Roth}
\affil[1]{Dept. of Mathematics and Statistics,
Univ. of Vermont,
Burlington, VT}
\affil[2]{Mitre Corp., Bedford, MA}
\affil[3]{Department of Mathematics, Emory Univ., Atlanta, GA}

%\vspace*{1.0ex}
\begin{comment}
\begin{center}
\begin{minipage}[t]{.45\textwidth}
   \begin{center}
     Dan Archdeacon \\
     Dept. of Math. and Stat. \\
     University of Vermont \\
     Burlington, VT 05405 \ \ USA\\
   \end{center}
 \end{minipage}\hspace*{.1\textwidth}
 \begin{minipage}[t]{.45\textwidth}
   \begin{center}
     Tom Boothby \\
     Dept. of Math. \\
     Simon Fraser University \\
     Burnaby, BC \ V5A\ 1S6  Canada\\
     {\tt tboothby@sfu.ca}
   \end{center}
 \end{minipage}
\end{center}
\vspace*{1.0ex}

\begin{center}
     Jeff Dinitz \\
     Dept. of Math. and Stat. \\
     University of Vermont \\
     Burlington, VT 05405 \ \ USA\\
     {\tt jeff.dinitz@uvm.edu}
\end{center}

\vspace*{1.0ex}
\begin{center}
{\bf Draft: Not for Distribution}\\
%\vspace*{1.0ex}
\today
\end{center}
\end{comment}

\maketitle
%\vspace*{1.0ex}
%\center{\today}

\begin{abstract}  The generalized taxicab number $T(n,m,t)$ is equal to the smallest number that is the sum of $n$ positive $m$th powers in $t$ ways. This definition is inspired by Ramanujan's observation that $1729 = 1^3+ 12^3 =9^3 + 10^3 $ is the smallest number that is the sum of two cubes in two ways and thus $1729= T(2,3,2)$.  In this paper we  prove that for any given positive integers $m$ and $t$, there exists a number $s$ such $T(s+k,m,t) =T(s,m,t) +k$ for every $k \geq 0$.  The smallest such $s$ is termed  the seed for the generalized taxicab number.  Furthermore, we find explicit expressions for this seed number when the number of ways $t$ is 2 or 3 and present a conjecture for $t \geq 4$ ways.
\end{abstract}

%%%%%%%%%%%%%%%%%%%%%%%%%%%%%%%%%%%%%%%%%%%%%%
\section {Introduction}\label{introduction}
%%%%%%%%%%%%%%%%%%%%%%%%%%%%%%%%%%%%%%%%%%%%%%

Hardy relays the following story about visiting Ramanujan during his illness (see \cite{Hardy}, p. xxxv):

\begin{center}
\begin{minipage}{5in}
{  I remember once going to see him when he was lying ill at Putney. I had ridden in taxi cab number 1729 and remarked that the number seemed to me rather a dull one, and that I hoped it was not an unfavorable omen. 'No,' he replied, 'it is a very interesting number; it is the smallest number expressible as the sum of two cubes in two different ways.' }
\end{minipage}  \end{center}

Indeed, $1729 = 1^{3} +12^{3}$ and $1729 = 9^{3} + 10^{3}$ and it is the smallest such number that is the sum of two cubes in two different ways.  In honor of the Ramanujan--Hardy conversation, the smallest number expressible as the sum of two cubes in $t$ different ways is known as the $t^{th}$ taxicab number and is denoted $\Taxi(t)$ .  Therefore, with this notation,  $\Taxi(2) = 1729$.

There has been quite a bit of effort expended in finding these taxicab numbers.  The interested reader is referred to \cite {Boyer} for information about these numbers. That paper contains interesting information about the history of the problem as well as  a discussion about the techniques used to find certain values of  $\Taxi(t)$. Further information about taxicab numbers and their variants can also be found at \cite{Boyer1,Carr,Meyrignac,Peterson,Silverman}.  Basically, $\Taxi(t)$ is known for $2\leq t\leq 6$ and upper bounds for $\Taxi(t)$ have been given for $7\leq t \leq 22$.

In this paper we will generalize the definition of taxicab numbers, as there is really nothing special about using exactly two cubes (except for historical reasons).  We will be concerned with finding the smallest number that is the sum of $n$ positive $m^{th}$ powers in at least $t$ ways.

Let $T(n,m,t) $ denote the least number that is the sum of $n$  positive $m^{th}$ powers in at least $t$ ways provided such a number exists.\footnote{It should be noted that in Wikipedia, the  generalized taxicab number $\Taxi(k, j, n)$ is the smallest number which can be expressed as the sum of $j$ $k^{th}$ positive powers in $n$ different ways, however since there have been no published papers with this notation we will use the notation in  the definition given above.}
So as noted above, $T(2,3,2)= 1729$ and in general $T(2,3,t)$ is the taxicab number $\Taxi(t)$.   It is also easy to verify that \\

\begin{comment}
$
\begin{array}{lll}
T(2,2,2)= 50& = 5^2+5^2& =7^2+1, \\

T(3,2,2)= 27&= 5^2+1+1 &=3^2+3^2+3^2, \\

T(4,2,2)= 28&= 5^2+1+1+1 &=3^2+3^2+3^2+1,   \ and\\

T(5,2,2)= 20& = 4^2+1+1+1+1 &= 2^2+2^2+2^2+2^2+2^2.\\
\end{array}$ 
\end{comment}

$T(2,2,2)= 50 = 5^2+5^2 =7^2+1, $ 

$T(3,2,2)= 27= 5^2+1+1 =3^2+3^2+3^2, $

$T(4,2,2)= 28= 5^2+1+1+1 =3^2+3^2+3^2+1,  $ \ and

$T(5,2,2)= 20 = 4^2+1+1+1+1 = 2^2+2^2+2^2+2^2+2^2.$\\

\noindent
Note that  by adding 1 to the two sums in the last example we  obtain
$$T(6,2,2)\leq  21 = 4^2+1+1+1+1+1 = 2^2+2^2+2^2+2^2+2^2+1.$$ Doing this again we obtain
$$T(7,2,2)\leq  22 = 4^2+1+1+1+1+1+1 = 2^2+2^2+2^2+2^2+2^2+1+1.$$ 

\noindent
In fact it is indeed true that $T(6,2,2)= 21$ and $T(7,2,2)= 22$.  \\

Considering cubes now, it is straightforward to verify that \medskip

$T(2,3,2)= 1729 = 12^3 +1 = 10^3+9^3,  $

$T(3,3,2)=  251= 6^3+3^3+2^3  = 5^3+5^3+1, $

$T(4,3,2)=  219 = 6^3 + 1+1+1= 4^3+4^3+4^3+3^3,$

$T(5,3,2)=  157 = 5^3+2^3+2^3+2^3+2^3 = 4^3+4^3+3^3+1+1 $

$T(6,3,2)=  158 = 5^3+2^3+2^3+2^3+2^3 +1= 4^3+4^3+3^3+1+1 +1$

$T(7,3,2)= 131= 5^3 + 1+1+1+1+1+1 = 4^3+3^3+2^3+2^3+2^3+2^3+2^3$

$T(8,3,2)=  132= 5^3 + 1+1+1+1+1+1 +1= 4^3+3^3+2^3+2^3+2^3+2^3+2^3+1$

$T(9,3,2)=  72 = 4^3 + 1+1+1+1+1+1 +1+1 = 2^3+2^3+2^3+2^3+ 2^3+2^3+2^3+2^3+2^3$

$T(10,3,2)=  73 = 4^3 + 1+1+1+1+1+1 +1+1 +1= 2^3+2^3+2^3+2^3+ 2^3+2^3+2^3+2^3+2^3+1 $\\

\noindent
One also can check that $T(11,3,2)=  74$ and that this solution comes about by adding 1 to both of the (equal) sums in the case of  $T(10,3,2)$.

From the examples above it seems plausible that there exists a number, say $s_0$, such that $T(n+1, m,t) = T(n, m,t) +1$  for all $n \geq s_0$  or equivalently, that $T(s_0+k, m,t) = T(s_0, m,t) +k$  for all $k \geq 0$ . This motivates the following definition.
\\

\noindent {\bf Definition.} If $s_0$ is the smallest positive integer such that $T(n+1, m,t) = T(n, m,t) +1$ for all $n \geq s_0$, then we call 
$s_0$ the {\em seed number} for $m^{th}$ powers in $t$ ways and denote this number by $S(m,t) = s_0$.   We also  call $T(S(m,t), m,t) $ the {\em seed value} of $m^{th}$ powers in $t$ ways and denote it by $V(m,t)$.
\medskip

%So, from the discussion and examples above we have that $S(2,2) = 5$ and $V(2,2) = 20$, and also that $S(3,2)= 9$ and $V(3,2) = 72$.

This paper proceeds as follows.
In Section  \ref{sect2} we will show  that for every $m$ and $t$ there exists a seed number.
In Section \ref{sect3} we will give an explicit value for the seed of the sum of $m^{th}$ powers in 2 ways.  We will prove there that  
the seed number for squares in 2 ways is indeed 5 and the seed value is 20.  In our notation, this says  $S(2,2) = 5$ and $V(2,2) = 20$ and hence $T(5+k, 2,t) = 20+k$ for all $k \geq 0$.  
In Section  \ref{sect4} we will give an explicit value of $V(m,3)$. In Section 5, we end with a general theorem and a conjecture about the seed for $m^{th}$ powers in $t$ ways for all $t \geq 2$.

\section {Seeds exist}\label{sect2}
  In this section we will prove that for every $m$ and $t$ there exists a seed number.  We will first show that for every $m$ and $t$ there exist some positive integer $n$  and some value $v$ such that $ v$ is the sum of $n$ $m^{th}$ powers in $t$ ways.  From that we will then show that there is a least such $n$ and hence there will exist a seed number (and a seed value).

\begin{lemma}\label{lemma1}For all positive integers $m,t \geq 1$, there exist positive integers $n$ and $v$  such that $ v$ is the sum of $n$ $m^{th}$ powers in $t$ ways. \end{lemma}

\proof
If $m=1$ or $t=1$, the result is obvious, so assume that $m,t >1$. 

 We will give a direct construction of $t$ different sums of $m^{th}$ powers, each sum  having the same number of terms.  The first sum will be $t^m + t^m + \ldots +t^m$ for a suitable number of terms.  For each $1\leq i\leq t-1$ we will construct the sum $S_i$ from the terms $(t-i)^m$ and $(t+i)^m$.  So each $S_i =   \underbrace {(t-i)^m + \ldots +(t-i)^m}_{x_i} + \underbrace{(t+i)^m + \ldots +(t+i)^m}_{y_i} $ for suitable $x_i$ and $y_i$.

  In order to define the sums, we first need to define some values. For each $1\leq i\leq t-1$ define 

\begin{tabular}{lll}
$a_i = t^m-(t-i)^m$ ;& $b_i = (t+i)^m-t^m$; & $\l_i = \mbox{lcm}(a_i,b_i)$; \\ 
$\alpha_i = l_i/a_i$; & $\beta_i = l_i/b_i$; & $\gamma_i = \alpha_i+\beta_i$; \\
$n=\mbox{lcm}(\gamma_1,\gamma_2, \ldots , \gamma_{t-1})$; &  $\delta_i = n/\gamma_i.$\\
\end{tabular}
\medskip

\noindent
Using these values, define
$$S_0 =   \underbrace {t^m + \ldots +t^m}_{n}$$ 
and for each $1\leq i\leq t-1$ define
 $$S_i =   \underbrace {(t-i)^m + \ldots +(t-i)^m}_{\alpha_i\delta_i} + \underbrace{(t+i)^m + \ldots +(t+i)^m}_{\beta_i\delta_i} .$$

%To define the $t$ multisets we will use  ``exponential'' notation, i.e. we write $\{a_1^{m_1}, a_2^{m_2} \ldots ,a_s^{m_s} \}$ to denote the multiset that contains the symbol $a_i$ $m_i$ times (for $1\leq i \leq s$).  We define    $ P_0 = \{ t^m\} = \{\underbrace{ t,t \ldots ,t}_{n}\}$  and for  $1\leq i\leq t-1$ define $P_i = \{ (t-i)^ {\alpha_i\delta_i},(t+i)^ {\beta_i\delta_i}\}$. 

Obviously all of these sums are different (i.e. no two contain the same terms).  We must prove two things: first, that each of the sums $S_0, S_1, \ldots, S_{t-1}$ contain the same number of terms and second,  that $S_0= S_1 = \ldots = S_{t-1}$.

To show the first,  we note  that $S_0$ contains $n$ terms and for each $1\leq i\leq t-1$,   $S_i$ contains $\alpha_i\delta_i + \beta_i\delta_i$ terms.  But now $$\alpha_i\delta_i + \beta_i\delta_i = (\alpha_i+ \beta_i)\delta_i = \gamma_i\delta_i = n$$
 as desired.

Next we compute the sums.  Clearly $S_0= nt^m$.  For each $1\leq i\leq t-1$,  

$$\begin{array}{lll}
S_i &= & \alpha_i\delta_i(t-i)^m + \beta_i\delta_i(t+i)^m\\
&= & (\alpha_i(t^m - a_i) + \beta_i(t^m+b_i))\delta_i\\
&= &  ((\alpha_i+\beta_i) t^m + (\beta_ib_i-\alpha_ia_i))\delta_i \\
&= &   ((\alpha_i+\beta_i) t^m + (l_i-l_i))\delta_i \\
&= &    \gamma_i t^m \ {n\over  \gamma_i}     \\
&= & nt^m 
\end{array}$$
Thus we conclude that $S_0= S_1 = \ldots = S_{t-1}= nt^m$. This completes the proof.
 \qed

\medskip

By adding a 1  to  each  multiset of terms in the above sums, we see that  $T(n+1,m,t) \leq T(n,m,t)$+1. 
As an immediate consequence of this observation and Lemma \ref{lemma1} we can conclude that  $T(n,m,t)$ exists  for every $m,t \geq 1 $ when   $n$ is large enough.

 In the next theorem we will prove that there is a seed number for every $m$ and $t$.
%% i.e. we will prove that for every $m,t >0$ there is a number $ s_0=S(m,t) $ such that  $T(s_0+k,m,t)= T(s_0,m,t) + k$ for every $k\geq 0$.

\begin{theorem} For every $m,t \geq 1$, there is a smallest number $s_0$ such that   $T(s_0+k,m,t)= T(s_0,m,t) + k$  for every $k\geq 0$. Hence for every $m,t \geq 1$ the seed number $s_0 = S(m,t)$ exists.
\end{theorem}

\proof Fix $m$ and $t$. 
Note that for some $n$ large enough, $T(n,m,t)$ exists by Lemma \ref{lemma1}.  Now let $n_0$ be the smallest integer for which $T(n_0,m,t)$ exists.
 To shorten notation let $T '(n) = T(n,m,t)$.  Notice that for any $n$,  $T'(n) \geq \underbrace{1^m+1^m +\ldots + 1^m}_n = n$.  So a (very naive) lower bound for   $T'(n)$ is $n$.  The {\em gap} between the value of  $T'(n_0)$  and the naive lower bound on $T'(n_0)$ is  $g = T'(n_0)-n_0 \geq 0$.  We also have that $T'(n_0+k) \leq T'(n_0)+k$ for every $k \geq 0$.  Thus $n_0 + k \leq T'(n_0+k) \leq T'(n_0)+k$.

We say that the function $T'$ {\em drops} at $n$ if $T'(n+1) < T'(n)+1$ (note that $n$ is the {\em location} of the ``drop'' and not the amount $T '$ drops).  Let $D = \{n \geq n_0 \ | \ T'(n+1) < T'(n)+1 \}$ be the set of all drops of $T'$. We claim that $D$ is a finite set  and in fact we will show that $|D| \leq g$.  Let $D= \{n_1, n_2, \ldots \}$. (We should note that possibly $D=\emptyset$, in which case  $S(m,t)$ is just $n_0$.)  If $D$ is nonempty, say that $n_k = n_0 + x_k$, then we see that since $T'$ drops by at least one at each $n_i$ we have that 
$$ T'(n_k+1) = T'(n_0+x_k+1) \leq T'(n_0) + x_k+1 - k.$$

So since $n_0+x_k+1 \leq T'(n_0+x_k+1)$, we have that $n_0 \leq T'(n_0)-k$ and hence that  $k\leq T'(n_0)-n_0 = g$.  This implies that $|D| \leq g$ and hence there are at most $g$ drops in the function $T'$.  So if $D = \{ n_1, n_2, \ldots n_i\}$ is the set of all drops, then  since $n_i$ is the last drop in the function $T'$, we have that $T'(n_i+1+k) = T'(n_i+1)+k$ for all $k \geq 0$ and hence $n_i+1$ is the seed number $S(m,t)$.  \qed

As an example of the above theorem we consider the values of  $T(2,3,2), T(3,3,2) ,\ldots ,$ $T(10,3,2)$ given in Section \ref{introduction}.  Notice that (from Ramanujan) $n_0=2$.  We also see that
$n_1 =2, n_2=3,n_3=4, n_4=6, $ and $n_5=8$.  We will prove in Section \ref{sect4} that indeed 
$D= \{2,3,4,6,8\}$.  Thus we will conclude that $S(3,2)=9$ and hence that $T(9+k,3,2) = T(9,3,2)+k = 72+k$ for all $k \geq 0$.

\section{Seeds for two  ways}\label{sect3}

In this section we will give the explicit value for the seed numbers for two  ways.  We assume that all variables are positive integers except where noted. We begin with three easy lemmas that hold for any number of sums.  The first  lemma says that if all the sums share a common term, then that term must be equal to 1.

\begin{lemma} \label{lemma3.1} If $x = T(n,m,t)$ and $x= \sum_{i=1}^n a_i^m =\sum_{i=1}^n b_i^m =\ldots = \sum_{i=1}^n t_i^m$, and if $a_{i_1}= b_{i_2} = \ldots = t_{i_t}$ for some choice of $i's$, then $a_{i_1}= b_{i_2} = \ldots = t_{i_t}= 1$.
\end{lemma}

\proof  If not, then replace each of $a_{i_1}, b_{i_2},\ldots ,t_{i_t}$ with a 1 and get a contradiction to $x = T(n,m,t)$ since the new sums of the $m^{th}$ powers will still all be the same and will be less than before, a contradiction.  \qed

The next lemma says that the $t$ sums adding to the seed value can't all have a 1 as a term.

\begin{lemma} \label{lemma3.2} 
If $x = V(m,t)$ and $x= \sum_{i=1}^n a_i^m =\sum_{i=1}^n b_i^m =\ldots = \sum_{i=1}^n t_i^m$, and if $a_{i_1}= b_{i_2} = \ldots = t_{i_t}$ for some choice of $i's$, then $a_{i_1}= b_{i_2} = \ldots = t_{i_t}\neq 1$.
\end{lemma}

\proof If each sum has a 1 as a term, then by simply deleting the 1 in each of these sums we would obtain a smaller seed value, a contradiction.  \qed

In the next lemma we show that  any seed value for $m^{th}$ powers in $t$ ways must always be greater or equal to $n2^m$ where $n$ is the seed number.  This will essentially say that  we can always assume that one of the sums is $ \underbrace{2^m+2^m + \ldots +2^m}_n $.  This fact will be of fundamental importance in finding seeds for 2 and 3 ways.

\begin{lemma} \label{powersof2}  If $V(m,t) = \sum_{i=1}^{s_0}a_i^m = \sum_{i=1}^{s_0}b_i^m =\ldots =\sum_{i=1}^{s_0}t_i^m$ is the seed value for $m^{th}$ powers in $t$ ways, then $V(m,t)\geq s_02^m$.  Further, if  $n$ is any number which provides a solution to 
$ \sum_{i=1}^na_i^m = \sum_{i=1}^nb_i^m =\ldots =\sum_{i=1}^nt_i^m= n 2^m$, then $V(m,t)= s_02^m$, where $s_0$ is the smallest such $n$ (and hence $s_0= S(m,t)$).

\end{lemma}

\proof  Let $V(m,t) = \sum_{i=1}^{s_0}a_i^m = \sum_{i=1}^{s_0}b_i^m =\ldots =\sum_{i=1}^{s_0}t_i^m$ be the seed value for $m^{th}$ powers in $t$ ways and  assume that  each sum is written in nonincreasing order.  
Now if $a_n =1$, then from Lemma \ref{lemma3.2} we have without loss of generality that $b_n \neq 1$. Thus for all $1\leq i \leq n$ it must be that $b_i \geq 2$.  Hence in this case we have that $V(m,t)\geq {s_0}2^m$.  If $a_n \geq 2$, then since $a_i \geq a_{i+1}$ for all $1\leq i \leq n-1$, then clearly $V(m,t)\geq {s_0}2^m$.
So we have that $V(m,2) \geq {s_0}2^m$.   

The second part of this lemma now follows immediately. \qed
\bigskip

 Lemma \ref{powersof2} says that the seed value is  the sum of $s_0$ $2^m$'s (where $s_0$ will be equal to the seed number $S(m,t)$).  So we are interested in this value for the sum.  In the next two lemmas we consider two different sums that are equal to the sum of $2^m$'s. The verification of the first is straightforward.

\begin{lemma} \label{4and1}  If $\alpha 4^m+(n-\alpha) = n2^m$, then $n=(2^m+1)\alpha$ and hence $n\geq(2^m+1)$.\end{lemma}

A comment is in order concerning  Lemma \ref{4and1}.  This lemma deals with the case when two sums are equal and one of the sums is all 2's and the other is 4's and 1's.  It says that if 
$$ \underbrace {4^m + 4^m + \ldots +4^m}_\alpha   + \underbrace{1+ 1+ \ldots +1}_{n-\alpha} = \underbrace{2^m + 2^m + \ldots +2^m}_n,$$
then $n\geq(2^m+1)$.  The next lemma deals with the case when two sums are equal and one of the sums is all 2's and the other is 3's and 1's. 

\begin{lemma} \label{3and1}  If $\alpha 3^m+(n-\alpha) = n2^m$ and if $d=\gcd(3^m-2^m,2^m -1)$, then  ${3^m-1\over d}\ | \ n$ and hence  $n \geq {3^m-1\over d}$. \end{lemma}

\proof  Assume that $\alpha 3^m+(n-\alpha) = n2^m$ and that $d=\gcd(3^m-2^m,2^m -1)$, then $$\alpha({3^m-1\over d}) = n({2^m-1\over d}).$$
Now since  $d=\gcd(3^m-2^m,2^m -1) =  \gcd(3^m-1,2^m -1)$, then $1= \gcd({3^m-1 \over d},{2^m -1\over d})$ and thus we have that ${2^m-1\over d}\ | \ \alpha.$  So $\alpha({d\over 2^m-1})$ is an integer.  Now since $\alpha({d\over 2^m-1}) \ ({3^m-1\over d}) = n$ it follows that $({3^m-1\over d}) \ |\ n$ and hence that $n \geq {3^m-1\over d}$.  \qed

The application of this lemma will be similar to that of Lemma \ref{4and1}.  In this case we have the situation where
$$ \underbrace {3^m + 3^m + \ldots +3^m}_\alpha   + \underbrace{1+ 1+ \ldots +1}_{n-\alpha} = \underbrace{2^m + 2^m + \ldots +2^m}_n.$$
So here we will have that $n= \alpha({d\over 2^m-1}) \ ({3^m-1\over d})$ and our main application will be that in this case  $n \geq {3^m-1\over d}$.
We now obtain our characterization of the seed number and the seed value for sums in two ways.

\begin{theorem}\label{2ways} Let $d= \gcd(3^m-2^m,2^m-1)$. 
The seed number $S(m,2) =\min({3^m-1\over d},2^m+1)=s_0$ and the seed value $V(m,2)=s_02^m$.   Hence $T(s_0+j,m,2) =s_02^m +j$   for every $j\geq 0$.
\end{theorem}

\proof   Consider the  two equations
$$4^m  + \underbrace{1+  \ldots +1}_{2^m}  = \underbrace{2^m +  \ldots +2^m}_{2^m+1}  = (2^m+1)2^m   \ \ \ \ \ \ (1)$$ and
$$\underbrace{3^m +\ldots + 3^m}_{2^m-1\over d} + \underbrace{1+ \ldots +1}_{3^m - 2^m\over d }  = \underbrace{2^m +  \ldots +2^m}_{3^m-1\over d} =({3^m-1\over d}) 2^m.   \ \ \ \ \ (2)$$

In view of Equations (1) and (2) and Lemma \ref{powersof2} if $V(m,t) = a_1^m+a_2^m+ \ldots +a_n^m = b_1^m+b_2^m+ \ldots +b_n^m$ is the seed value for $m^{th}$ powers in 2 ways, then we can assume that $ b_i=2$ for all $i$. We next show that $a_1 \leq 4$.

Assume that $V(m,t) = a_1^m+a_2^m+ \ldots +a_n^m = 2^m+2^m+ \ldots +2^m$ is the seed value for $m^{th}$ powers in 2 ways with $a_i \geq a_{i+1}$  for all $1\leq i \leq n-1$. 
Assume that $a_1 \geq 5$.  Clearly, if $n \geq 2^m+1$, then in view of Equation (1) above this is a contradiction (since $a_1^m+a_2^m+ \ldots +a_n^m > 4+1+ \ldots +1$).  Assume $n <2^m+1$, then
$$a_1^m+a_2^m+ \ldots +a_n^m + \underbrace{1 + 1+\ldots +1}_{2^m+1-n}  > 4^m  + \underbrace{1+  \ldots +1}_{2^m}$$ which is again a contradiction to the assumption that $V(m,t) = a_1^m+a_2^m+ \ldots +a_n^m $.  Assuming that $a_1=4$ and $a_2 >1$ yields a similar contradiction.

So either $\{a_1,a_2, \ldots a_n\} = \{1,4\}$ or $\{a_1,a_2, \ldots a_n\} = \{1,3\} $ since by Lemma \ref{lemma3.1} no $a_i$ can equal any $b_i=2$.  In the first case we obtain Equation (1) since no smaller sum can have only $4^m$'s and 1's as its terms. In the second case we can assume that 
$s 3^m+(n-s) = n2^m$  for some $s$.  
From Lemma \ref{3and1}, the minimum value of $n$ is $ {3^m-1\over d}$, which leads to Equation (2).  The seed number will therefore be the minimum length of the sums in either Equation (1) or Equation (2). Thus
the minimum of $2^m+1$ and $3^m-1\over d$ will be the seed number $S(m,2)$. \qed

In the following table we compute seeds for  $m^{th}$ powers in 2 ways for $m \leq 20$.

\begin{example}  In the following table we give explicit values from Theorem \ref{2ways}.  Remember that $S(m,2)$ is the number of terms in the seed, while $V(m,2)$ is the exact value of the seed.
\end{example}

\begin{center}
\begin{tabular} {|lccc|}\hline
$m$&$d$&$S(m,2)$&$V(m,2)$\\ \hline
1 &1 &2 & 4 \\
2 &1 &5 & 20 \\ 
3 &1 &9 & 72 \\ 
4 &5 &16 & 256 \\ 
5 &1 &33 & 1056 \\ 
6 &7 &65 & 4160 \\ 
7 &1 &129 & 16512 \\ 
8 &5 &257 & 65792 \\ 
9 &1 &513 & 262656 \\ 
10 &11 &1025 & 1049600 \\ 
11 &23 &2049 & 4196352 \\ 
12 &455 &1168 & 4784128 \\ 
13 &1 &8193 & 67117056 \\ 
14 &1 &16385 & 268451840 \\ 
15 &1 &32769 & 1073774592 \\ 
16 &85 &65537 & 4295032832 \\ 
17 &1 &131073 & 17180000256 \\ 
18 &133 &262145 & 68719738880 \\ 
19 &1 &524289 & 274878431232 \\ 
20 &275 &1048577 & 1099512676352 \\ \hline
\end{tabular} 
\end{center}

The interested reader may note that $S(m,2) = 2^m+1$ in every case above except when $m=1,4,12$.  This says that $2^m+1 \leq (3^m-1)/d$ for every $m \leq 20$ with $m \neq 1,4,12$.    We computed values of $2^m+1$ and  $(3^m-1)/d$ for all $m \leq 200,000$ and found that $2^m+1 \leq (3^m-1)/d$ for all $m$ in that range except for $m=1,4,12$ and 36.  We do not conjecture that this is true for all $m >36$, it sure appears to be true, however.
%\marginpar{want to make this conjecture?}

\section{Seeds for three  ways}\label{sect4}

In this section we will give an explicit value for $V(m,3)$, the seed value for the smallest number that can be written as the sum of $m^{th}$ powers in 3 ways.
We first need a preliminary lemma that says that no term in a sum that is a seed value (for  powers $m \geq 4$) can exceed the number 4.

\begin{lemma}\label{no5} The seed value $V(m,3) = \sum_{i=1}^s a_i^m =\sum_{i=1}^s b_i^m = s2^m$ where $s = S(m,3)$ is the seed number. If $m\geq 4$, then $a_i,b_i \leq 4$ for all $1\leq i\leq s$.\end{lemma}

\proof  We first note following equation:
{\footnotesize $$ \underbrace{4^m+ 1+ \ldots +1}_{2^m+1} +  \underbrace{4^m+ 1+ \ldots +1}_{2^m+1} = \underbrace{4^m+ 1+ \ldots +1}_{2^m+1} +  \underbrace {2^m+  \ldots +2^m}_{2^m+1}= \underbrace {2^m+  \ldots +2^m}_{2^m+1}+\underbrace {2^m+  \ldots +2^m}_{2^m+1} . \mbox{ (3)}$$}
From this equation and Lemma \ref{powersof2} we have $V(m,3) = s 2^m$ for  $s = S(m,3)$. We also see from this equation 
that the taxicab number $T(2(2^m+1)+j, m,3) \leq 2(2^m+1)2^m+j$ for all $j \geq 0$ and that the seed number $S(m,3) \leq 2(2^m+1)$. So in particular, when $j=0$ we have that $ T(2(2^m+1),m,3) \leq 2(2^m+1)2^m$.

Now, assume $V(m,3) = \sum_{i=1}^s a_i^m =\sum_{i=1}^s b_i^m = s2^m$ where $s = S(m,3)$ is the seed number. Then $s \leq 2(2^m+1)$.  Assume that  $a_1 \geq 5$.  Since $a_i \geq 1$ for all $i >1$, when extending the sums to   have $2(2^m+2)$ terms by adding sufficinetly many 1's, we get that
$$5^m + 2(2^m+1)-1 \leq \sum_{i=1}^s a_i^m +   (2(2^m+1)-s) \leq T(2(2^m+1),m,3)$$
and hence 
$$5^m + 2(2^m+1)-1 \leq T(2(2^m+1),m,3) \leq 2(2^m+1)2^m.$$
Thus
\begin{center}$
\begin{array}{rcll}
5^m + 2(2^m+1)-1 &\leq & 2(2^m+1)2^m&\\
5^m + 2(2^m+1)-1 &\leq & (2^m+1)2^{m+1}&\\
5^m &\leq& (2^m+1)(2^{m+1}-2) + 1& \\ 
5^m &\leq& 2^{2m+1}-1&\\
5^m &\leq& 2\times 4^m-1.&\\
\end{array}$
\end{center}

This last inequality implies that $m = 1,2,$ or 3, but   by hypothesis  $m \geq 4$, so we obtain a contradiction.  Hence 
 $a_i \leq 4$ (similarly $b_i\leq 4$) for all $1\leq i\leq s$. \ \ \qed

\bigskip
We are now in position 
to obtain our characterization of the seed number and the seed value for sums in three ways.  We begin with the small values of $m$.

\begin{theorem}\label{small3way}(a)  $S(1,3) = 3$ and the seed value $V(1,3)=3 \times 2^1 =6$, 
(b) $S(2,3) = 8$ and the seed value $V(2,3)=8 \times2^2 =32$, (c) $S(3,3) = 18$ and the seed value $V(3,3)=18\times 2^3 =144$.
  \end{theorem}

\proof The sums are given below. It is straightforward to check that they are minimal.\\

\noindent
(a)\ $6=4+1+1 = 3+2+1= 2+2+2$ \\ \\
(b)\ $32 = 4^2+ 2^2+2^2+2^2 + 1+1+1+1 = 3^2+3^2+3^2 + 1+1+1+1+1 = \underbrace{2^2 +\ldots + 2^2}_8  $  \\
(c)\ $144= 4^3+4^3 + \underbrace{1 + \ldots +1}_{16} = 4^3+ \underbrace{1 + \ldots +1}_8 + \underbrace{2^3 +\ldots + 2^3}_9 = \underbrace{2^3 +\ldots + 2^3}_{18}$

\begin{theorem}\label{3ways}  
Assume that $m \geq 4$ and let $d= \gcd(3^m-2^m,2^m-1)$.  Also let $l_3 = {3^m-1\over d}$ and $l_4 =2^m+1$. 
Given the four values $l_3,l_4,2l_3,2l_4$, the second smallest of these values is the seed number $S(m,3)$ and the seed value 
$V(m,3) = S(m,3) \times 2^m.$ 

\begin{comment}
Then 
\begin{enumerate}
  \item if\ $2l_4<l_3$, then the seed number $S(m,3) = 2l_4 $ and the seed value $V(m,3)=2l_42^m $.
  \item if\ $l_3 \leq l_4 \leq 2 l_3$, then the seed number $S(m,3) = l_4 $ and the seed value $V(m,3)= l_42^m$,    
  \item if\ $l_4 <l_3 \leq 2 l_4$, then the seed number $S(m,3) = l_3 $ and the seed value $V(m,3)= l_32^m$, 
  \item if\ $2l_3 <l_4$, then the seed number $S(m,3) = 2l_3 $ and the seed value $V(m,3)=2l_32^m $, 
\end{enumerate}
\end{comment}

\end{theorem}

\proof
Considering Equation (3) in the proof of Lemma \ref{no5}, in all cases the seed number $S(m,3) \leq 2l_4 $. Also, we can assume that the  seed value $V(m,3) = \sum_{i=1}^s a_i^m =\sum_{i=1}^s b_i^m = s2^m$ where $s = S(m,3)$ is the seed number and (from Lemma \ref{no5}) that $a_i, b_i \leq 4$ for all $i$.  Let $A = \{a_1,a_2, \ldots a_s\}= \{1^{\alpha_1},2^{\alpha_2},3^{\alpha_3},4^{\alpha_4}\} $ be the multiset containing all the terms in the sum $\sum_{i=1}^s a_i^m$ (so $A$ contains the term $i^m$ exactly $\alpha_i$ times for $1\leq i\leq 4$), and let $B = \{b_1,b_2, \ldots b_s\}= \{1^{\beta_1},2^{\beta_2},3^{\beta_3},4^{\beta_4}\} $ be the multiset containing  the terms in the sum $\sum_{i=1}^s b_i^m$.  From Lemma \ref{lemma3.1} we can assume without loss of generality that $\beta_2=0$. 

\medskip\noindent
{\em Case 1.)} \ \ If $2l_4$ is the second smallest value, then  $2l_4\leq l_3$. Assume that  $s=S(m,3) < 2l_4$.   We see first that $0\leq \alpha_4,\beta_4 \leq 1$, 
since if (say) $\alpha_1 \geq 2$, then $\sum_{i=1}^s a_i^m \geq 4^m+4^m +(s-2)1^m$ and so $\sum_{i=1}^s a_i^m + (2l_4-s) \geq 4^m+4^m +(2l_4-2)$ which (because of Equation (3)) says that 
$\sum_{i=1}^s a_i^m$ can  not be a seed value unless $\alpha_3=\alpha_2=0$ in which case we are led to one of the sums in Equation (3).  However since we assumed that $s < 2l_4$ this is a contradiction.

Now, if $\beta_4 = 0$ we obtain the equation $\beta_3 3^m +\beta_1 1^m = s2^m$.  By Lemma \ref{3and1} we thus have that $s \geq l_3$.  So $s \geq l_3> 2l_4>s$ a contradiction.  If $\alpha_4 = 0$, then by subtracting $\alpha_2$ $2^m$'s from each side of the equation $\sum_{i=1}^s a_i^m = s2^m$ we obtain a similar contradiction.  Hence we can assume that $\alpha_4= \beta_4 = 1$.

So we have that $$4^m + \alpha_3 3^m + \alpha_2 2^m + \alpha_1 1^m =  4^m + \beta_3 3^m +  \beta_1 1^m$$
subtracting $4^m$ from both sides yields
$$ \alpha_3 3^m + \alpha_2 2^m + \alpha_1 1^m =   \beta_3 3^m +  \beta_1 1^m.$$
This implies that 
$$\alpha_22^m = (\beta_3-\alpha_3)3^m + (\beta_1 -\alpha_1)$$
Now since $(\beta_3-\alpha_3)+  (\beta_1 -\alpha_1) = \alpha_2$ and since $\alpha_2 >0$ (else the equation is degenerate), then by Lemma \ref{3and1} we have that $\alpha_2 \geq l_3$, a clear contradiction.  

So in this case we have that $s \geq 2 l_4$.  Equation (3) then proves that indeed in this case that $s=S(m,3) = 2l_4$ and hence the seed value $V(m,3) = 2l_4 2^m$.

\medskip\noindent
{\em Case 2.)} If $l_4$ is the second smallest value, then    $l_3 \leq l_4 \leq 2 l_3$.
First consider the following equation:
{ $$ \underbrace{4^m+ 1+ \ldots +1}_{l_4}  = \underbrace{3^m+ \ldots +3^m}_{(2^m-1)/d)} + \underbrace{1^m+ \ldots +1^m}_{(3^m-2^m)/d)}+ \underbrace {2^m+  \ldots +2^m}_{l_4-l_3}= \underbrace {2^m+  \ldots +2^m}_{\l_4} .\ \ \  \mbox{ (4)}$$}
From this equation we see that in this case $s=S(m,3) \leq l_4$.  Assume  $s=S(m,3) < l_4$.

If $\alpha_4\geq1$ (or $\beta_4\geq1$), then  $\sum_{i=1}^s a_i^m \geq 4^m +(s-1)1^m$ and so $\sum_{i=1}^s a_i^m + (l_4-s) \geq 4^m +(l_4-1)$ which (because of Equation (4)) says that 
$\sum_{i=1}^s a_i^m$ can  not be a seed value, unless $\alpha_3=\alpha_2=0$ in which case we are led to the first  sum in Equation (4). However, since we assumed that $s < l_4$ we see that this is a contradiction. So $\alpha_4=\beta_4=0.$

Hence  
 $$ \alpha_3 3^m + \alpha_2 2^m + \alpha_1 1^m =  \beta_3 3^m +  \beta_1 1^m = s2^m \ \ \ \ \ \  \ (5)$$
subtracting $\alpha_3 3^m $ and  $\alpha_1 1^m$ from the first two sums yields
$$\alpha_2 2^m = (\beta_3-\alpha_3) 3^m +  (\beta_1-\alpha_1) 1^m$$
and so by Lemma \ref{3and1} we have that $\alpha_2 \geq l_3.$  Also, subtracting $\alpha_2 2^m$ from  the first and third
sums in Equation (5) we have 
$$ \alpha_3 3^m + \alpha_1 1^m = (s-\alpha_2) 2^m.$$
So, again by Lemma \ref{3and1} we have
 $s-\alpha_2 = \alpha_1+\alpha_3 \geq l_3$.  Thus $s= \alpha_1+\alpha_2+\alpha_3 \geq 2l_3$ a clear contradiction to our assumption that $s <l_4 \leq 2 l_3$.

So in this case we have $s \geq l_4$.  Equation (4) then proves that in this case that $s=S(m,3) = l_4$ and hence the seed value $V(m,3) = l_4 2^m$.

\medskip\noindent
{\em Case 3.)}      If $l_3$ is the second smallest value, then we have   $l_4 \leq l_3 \leq 2 l_4$.  First note the following equation:
{ $$ \underbrace{4^m+ 1+ \ldots +1}_{l_4} +\underbrace{2^m + \ldots +2^m}_{l_3-l_4}= \underbrace{3^m+ \ldots +3^m}_{(2^m-1)/d)} + \underbrace{1^m+ \ldots +1^m}_{(3^m-2^m)/d)}= \underbrace {2^m+  \ldots +2^m}_{\l_3} .\ \ \  \mbox{ (6)}$$}
From this we see that $s=S(m,3) \leq l_3$.  Assume that $s=S(m,3) < l_3$. So
similar to Case 1, the fact that $s < 2l_4$ implies that $\alpha_4,\beta_4 \leq 1$. We show this by considering $T(2l_4,m,3)$. This value is equal to $  \sum_{i=1}^s a_i^m + (2l_4 - s)$ since $\sum_{i=1}^s a_i^m$ is the seed value.  Note $\sum_{i=1}^s a_i^m + (2l_4 - s)  < 2l_4 2^m$. However, if $\alpha_4 \geq 2$, then $4^m+4^m + (l_4-2)1^m = 2l_4 2^m$, a contradiction.  

So we can assume that $0\leq\alpha_4,\beta_4 \leq 1$. We have that
$$ \alpha_4 4^m+ \alpha_3 3^m + \alpha_2 2^m + \alpha_1 1^m =  \beta_4 4^m+\beta_3 3^m +  \beta_1 1^m = s2^m.$$
If $\beta_4 = 0$, then we have $\beta_3 3^m +  \beta_1 1^m = s2^m $, but from Lemma \ref{3and1} this implies that $ s\geq l_3$, a contradiction our assumption that $s< l_3$.  So now we have that 
$$ \alpha_4 4^m+ \alpha_3 3^m + \alpha_2 2^m + \alpha_1 1^m = 4^m+ \beta_3 3^m +  \beta_1 1^m = s2^m. $$
If $\alpha_4 = 1$, then  by subtracting $4^m$ from the first two sums in the equation above, we obtain $\alpha_2 2^m =(\beta_3-\alpha_3)3^m +(\beta_1-\alpha_1)1^m$. But from Lemma \ref{3and1} we have  $\alpha_2 \geq l_3$, a contradiction. So $\alpha_4 = 0$. So now we have
$$  \alpha_3 3^m + \alpha_2 2^m + \alpha_1 1^m = 4^m+ \beta_3 3^m +  \beta_1 1^m = s2^m. $$
Finally, from   $  \alpha_3 3^m + \alpha_2 2^m + \alpha_1 1^m = s2^m $ subtract $\alpha_2 2^m$ from each side to obtain
$  \alpha_3 3^m +  \alpha_1 1^m = (s-\alpha_2)2^m $.  From Lemma \ref{3and1} we obtain $s-\alpha_2 \geq l_3$ a clear contradiction to our assumption that $s< l_3$.  

So we have shown that $s \geq l_3$.  Equation (6) then proves that in this case that $s=S(m,3) = l_3$ and hence the seed value $V(m,3) = l_3 2^m$.

\medskip\noindent
{\em Case 4.)} Finally,  if $2l_3$ is the second smallest value, then $2l_3 \leq l_4$.  Consider the following equation:
$$ \underbrace{3^m+ \ldots +3^m}_{(2^m-1)/d)} + \underbrace{1^m+ \ldots +1^m}_{(3^m-2^m)/d)}+ \underbrace{3^m+ \ldots +3^m}_{(2^m-1)/d)} + \underbrace{1^m+ \ldots +1^m}_{(3^m-2^m)/d)} = \ \ \ \ \ \ \ \  $$
 $$\underbrace{3^m+ \ldots +3^m}_{(2^m-1)/d)} + \underbrace{1^m+ \ldots +1^m}_{(3^m-2^m)/d)}+\underbrace {2^m+  \ldots +2^m}_{\l_3} = 
\underbrace {2^m+  \ldots +2^m}_{2\l_3}. \hspace {.5in} \mbox{(7)} $$
From this we see that $s=S(m,3) \leq 2l_3$.  Assume that $s=S(m,3) < 2l_3$.

Since $s<l_4$ (as in the proof of Case 2) we have that $\alpha_4=\beta_4 =0$.
So we have that
 $$ \alpha_3 3^m + \alpha_2 2^m + \alpha_1 1^m =  \beta_3 3^m +  \beta_1 1^m = s2^m. \ \ \ \ \ \  \ \ \ (8)$$
Subtracting $\alpha_3 3^m $ and  $\alpha_1 1^m$ from the first two sums yields
$$\alpha_2 2^m = (\beta_3-\alpha_3) 3^m +  (\beta_1-\alpha_1) 1^m$$
and so by Lemma \ref{3and1} we have that $\alpha_2 \geq l_3.$  Also, subtracting $\alpha_2 2^m$ from  the first and third
sums in Equation (8) we have 
$$ \alpha_3 3^m + \alpha_1 1^m = (s-\alpha_2) 2^m.$$
So, again by Lemma \ref{3and1} we have
 $s-\alpha_2 = \alpha_1+\alpha_3 \geq l_3$.  Thus $s= \alpha_1+\alpha_2+\alpha_3 \geq 2l_3$ a contradiction to our assumption that $s \leq 2 l_3$.

So in this case $s \geq 2l_3$.  Equation (7) then proves that in this case  $s=S(m,3) = 2l_3$ and hence the seed value $V(m,3) = l_4 2^m$.   \qed
\bigskip

As was done after the proof of Theorem \ref{2ways} we wish to compute the exact value of $S(m,3)$ using the results of Theorem \ref{3ways}.
We found that for every $1\leq m \leq 200,000$ with $m \neq 1,2,4,6,12,36$,  that $S(m,3) = 2l_4$ and hence that $V(m,3)= 2l_42^m.$ 
This  is Case 1 above and says that $2 l_4 \leq l_3$ for all powers $36< m \leq 200,000$.  Again we do not conjecture that $2 l_4 \leq l_3$ for all $m >36$, but certainly the evidence is very strong.

  When  $m=4$ we are in Case 2, so $l_3 \leq l_4 \leq 2 l_3$ and hence $S(4,3) = l_4= 2^4+1$.  We should note that $m=1$ also has the property that $l_3 \leq l_4 \leq 2 l_3$, and although it doesn't follow from the general proof, it is indeed true that $S(1,3) = l_4 = 3$ and so $V(1,3) = 3 \times 2^1=6$, since $6=1+2+3=4+1+1=2+2+2.$

When $m=6$ we are in Case 3 so  $l_4 \leq l_3 \leq 2 l_4$, and hence $S(6,3) = l_3 = 104$.
Also note that $m=2$ also has the property that $l_4 \leq l_3 \leq 2 l_4$, and  indeed $S(2,3) = l_3 = 8$ and so $V(1,3) = 8 \times 2^2=32$, since 
$32 = 4^2+ 2^2+2^2+2^2 + 1+1+1+1 = 3^2+3^2+3^2 + 1+1+1+1+1 = 2^2 + 2^2 + 2^2 + 2^2 + 2^2 + 2^2 + 2^2 + 2^2  $. As a bonus here we see that  
$32 = 5^2 + 7 \times 1^2$ and hence four different sums of 8 squares are equal to 32.  So we get  $S(2,4) = 8$ and hence $T(8,2,4)=
T(8,2,3) =32$.

Finally when $m=12$ or $m=36$, we have  $2 l_3 < l_4$ and thus both these values fall into case 4.

\section{More than 3 ways}

In this section we present  a general theorem and a conjecture about seeds.    Again  let $d= \gcd(3^m-2^m,2^m-1)$, $l_3 = {3^m-1\over d}$ and $l_4 =2^m+1$. Also let the sums
$$S_1= \underbrace{3^m +\ldots + 3^m}_{2^m-1\over d} + \underbrace{1+ \ldots +1}_{3^m - 2^m\over d }  \mbox{   and   }
S_2 =4^m  + \underbrace{1+  \ldots +1}_{2^m}  .$$

\begin{theorem}\label{generaltheorem}  
Given $t$, there exists a number $m_0$, such that if $m \geq m_0$, then the seed number $S(m,t)$ is bounded above by
the $t-1$st  smallest of the values $al_3+bl_4$ over all  $a, b \geq 0$.   The seed value 
$V(m,t) = S(m,t) \times 2^m.$ 
\end{theorem}

\proof  Let $n = \min(l_3,l_4)$ and let $m_0 = \max \{ m \ |\ 5^m< (t-1) n2^m\}$.   Let $n_0$ be the $t-1$st smallest of the values $al_3+bl_4$ over all  $a, b \geq 0$.  Finally define the sum $aS_1+bS_2+\overline{2^m}$ to be $a$ copies of $S_1$ added to $b$ copies of $S_2$ added to $n_0 - (al_3+bl_4)$ 
copies of $2^m$. Now it is clear that for all $al_3+bl_4 \leq n_0$, that the sum $aS_1+bS_2+\overline{2^m} = n_02^m$.   Thus the $t-1$ sums $aS_1+bS_2+\overline{2^m}$ with the sum of $m_0$ $2^m$'s are all equal, proving our upper bound. 
Furthermore, we note that if $m \geq m_0$, then   $k^m \geq (t-1) n2^m \geq n_02^m$ for all $k\geq 5$ and hence no sum with  $n_0$ terms and equal to  $n_02^m$ can contain any $k^m$ for $k\geq 5$. \ \qed
  
Note that in the proof of the previous theorem we did not need to prove that no sum could contain a $k^m$ for any $k \geq 5$ in order to obtain an upper bound.  We included that fact in order to add credence to our conjecture below.
Indeed we believe that the number presented in Theorem \ref{generaltheorem} is the actual seed number. We state this in the following conjecture.  One can see that both Theorem \ref{2ways} and Theorem \ref{3ways}  follow from this conjecture.

\begin{conjecture}\label{conjecture}  
Given $t$, there exists a number $m_0$, such that if $m \geq m_0$, then 
the $t-1$st  smallest of the values $al_3+bl_4$ over all  $a, b \geq 0$ is the seed number $S(m,t)$ and the seed value 
$V(m,t) = S(m,t) \times 2^m.$ 
\end{conjecture}

\section{Conclusion}  The generalized taxicab number $T(n,m,t)$ is equal to the smallest number that is the sum of $n$ $m^{th}$ powers in $t$ ways. This definition is inspired by Ramanujan's observation that $1729 = 1^3+ 12^3 =9^3 + 10^3 $ is the smallest number that is the sum of two cubes in two ways and thus $1729= T(2,3,2)$.  In this paper we first proved that for any given positive integers $m$ and $t$, there exist a seed for the generalized taxicab number, i.e. there exists a number $s=S(m,t)$ such $T(s+k,m,t) =T(s,m,t) +k$ for every $k \geq 0$.    We then found explicit expressions for this seed number when the number of ways $t$ is 2 or 3.  We ended  with a general theorem and conjecture about the seed number $S(m,t)$ for all $t$.

\bigskip \bigskip
\noindent {\bf Addendum:} Research for this paper was mostly undertaken while the authors were together in their first year of graduate school in mathematics at The Ohio State University in 1974.  This paper should have appeared shortly after that, but at least it is finally finished now. (It only took another 45 years).

%%%%%%%%%%%%%%%%%%%%%%%%%%%%%%%%%%%%%%%%%%%%%%
%%%%%%%%%%%%%%% THE REFERENCES %%%%%%%%%%%%%%%%%%
%%%%%%%%%%%%%%%%%%%%%%%%%%%%%%%%%%%%%%%%%%%%%%

------------------------------------------------------------------------------------------------------------------

(Concerned with sequence A011541 .)

------------------------------------------------------------------------------------------------------------------

%%%%%%%%%%%%%%%%%%%%%%%%%%%%%%%%%%%%%%%%%%%%%%
%%%%%%%%%%%%%%%%%%%%%%%%%%%%%%%%%%%%%%%%%%%%%%
%%%%%%%%%%%%%%%%%%%%%%%%%%%%%%%%%%%%%%%%%%%%%%

\end{document}